\documentclass[12pt,leqno]{article}
\usepackage[latin1]{inputenc}
\usepackage{amsmath,amsfonts,amssymb,amsthm,amscd}
\usepackage{calrsfs}
\usepackage{color}
\date{}

\newcommand{\Z}{\mathbb{Z}}

\newcommand{\A}{\mathbb{A}}

\newcommand{\C}{\mathbb{C}}

\newcommand{\G}{\mathbb{G}}

\newcommand{\R}{\mathbb{R}}

\newcommand{\Q}{\mathbb{Q}}

\newcommand{\Ss}{\mathfrak{S}}

\newcommand{\fn}{\mathfrak{n}}
\newcommand{\fm}{\mathfrak{m}}
\newcommand{\fu}{\mathfrak{u}}
\newcommand{\gog}{\mathfrak{g}}

\newcommand{\Hh}{\mathcal{H}}
\newcommand{\Aa}{\mathcal{A}}

\newcommand{\Ll}{\mathcal{L}}

\newcommand{\Oo}{\mathcal{O}}

\newcommand{\CSs}{\mathcal{S}}

\newcommand{\Res}{\mathrm{Res}}
\newcommand{\ind}{\mathrm{ind}}

\newcommand{\rg}{\rightarrow}
\newcommand{\lgr}{\longrightarrow} 
\newcommand{\Aut}{\mathrm{Aut}}  
\newcommand{\Gal}{\mathrm{Gal}}

\newcommand{\Ind}{\mathrm{Ind}}

\newcommand{\resp}{\mathrm{resp}}

\newcommand{\Alg}{{\mathrm{Alg}}}
\newcommand{\rk}{{\mathrm{rk}}}
\newcommand{\uU}{{\mathrm{U}}}
\newcommand{\St}{{\mathrm{St}}}
\newcommand{\GL}{{\mathrm{GL}}}
\newcommand{\SL}{{\mathrm{SL}}}
\newcommand{\SO}{{\mathrm{SO}}}
\newcommand{\Sp}{{\mathrm{Sp}}}
\newcommand{\usp}{{\mathrm{sp}}}
\newcommand{\uO}{{\mathrm{O}}}
\newcommand{\trace}{{\mathrm{trace}}}
\newcommand{\sgn}{{\mathrm{sgn}}}

\newcommand{\ba}{\backslash}

\newcommand{\Lie}{\mathrm{Lie}}

\newcommand{\IH}{IH}

\newcommand{\nequiv}{\equiv\kern-4.5mm/}

\def\adots{\mathinner{\mkern1mu\raise1pt\vbox{\kern7pt\hbox{.}}
\mkern2mu\raise4pt\hbox{.}
\mkern2mu\raise7pt\hbox{.}\mkern1mu}}

\newcommand{\intt}{\int\kern-4.5mm-}
\newcommand{\liminj}{\varinjlim}

\newtheorem{monlem}{{Lemma}}[section]
\newtheorem{madef}{{Definition}}[section]
\newtheorem{maprop}{{Proposition}}[section]
\newtheorem{maconj}{{Conjecture}}[section]
\newtheorem{montheo}{{Theorem}}[section]

\numberwithin{equation}{section}

\begin{document}

\title{On the central value of Rankin $L$-functions for self-dual algebraic representations of linear groups over totally real fields}

\author{Laurent Clozel and Arno Kret}

\maketitle

\section*{Introduction}

\textbf{1.} Let $F$ be a number field. Grothendieck has defined the category of pure motives over $F$. (See \cite{Dem}, \cite{Ma} for the definitions. In this section, which is aimed at motivating our results, we assume the optimistic properties of motives, defined using the relation of numerical equivalence of cycles.) 

Assume $E$ is another number field. There is a natural notion of motives with coefficients in $E$. If $\sigma:E\rg \C$ is an embedding, one conjectures the existence of an $L$-function $L(\sigma,M,s)$ having the usual properties, analytic continuation and functional equation. See \cite[\S~1]{Del}.

From a motive $M$ over $E$ we can deduce, for $E'\supset E$, a motive $M_{E'}=M \otimes E'$. For a suitable choice of $E'$, $M_{E'}$ will be absolutely irreducible. For $\sigma:E'\rg \C$, it is conjectured that $L(\sigma,M_{E'},s)$ is equal to the $L$-function $L(s,\pi)$ of a cuspidal (complex) representation of $\GL(N,\A_F)$  where $N$ is the dimension of $M_{E'}$. This conjecture was formulated  in the  70's, see Langlands \cite{Lan} and Serre \cite{Se}.

The class of cuspidal representations so obtained was specified by Borel \cite{Bo}, and more precisely in \cite{Cl} where the correspondence was precisely spelled out \cite[Ch.~4]{Cl}.

Now write simply $M$ for an absolutely irreducible (pure) motive with $E$-coefficients. In \cite{Del} Deligne has specified a set of critical integers and given a conjecture for the value $L(M,\sigma,n)$ for $n$ a critical integer. It should be, up to an element of $E^\times$, a \textit{period} given by the relation between the Betti and the de Rham cohomology of $M$. In particular, under this conjecture
\begin{equation}
L(M,\sigma,n) = 0 \Leftrightarrow L(M,\sigma',n)=0
\end{equation}
for two embeddings $\sigma,\sigma'$. (Gross further conjectured that the order of vanishing should be the same.)

\vspace{2mm}

\noindent\textbf{2.} Given the conjectural correspondence, this can be formulated unconditionally for the $L$-functions of algebraic cuspidal representations, as defined in \cite{Cl}. Consider  such a representation $\pi$ of $\GL(N,\A_F)$. From now on we assume $\pi$ self-dual and $N$ even. We also assume $\pi$ algebraic and regular in the sense of \cite{Cl}. Then the finite component $\pi_f$ of $\pi$ is defined over a number field $E$ \cite[Theorem 3.3]{Cl}\footnote{ In the statement of the Theorem, $\pi$ should be in $\Alg^0(n)$, the set of \textbf{cuspidal} algebraic representations.}. In particular, for $a \in \Aut(\C)$, there exists a representation $a(\pi)$ with finite part $a(\pi)_f = a(\pi_f)$.

For $N$ even, the correspondence between $L(\sigma,M,s)$ and $L(s,\pi)$ --- where $\pi$ should be associated to $(M,\sigma)$ --- is not trivial. Indeed, the usual multiplication by $q_v^{1/2}$ in the formation of the $L$--function intervenes (Hecke, cf. Tate \cite{Ta}). This is detailed in \S~1.4. It follows that the critical values for $L(s,\pi)$ are contained in $\Z+\frac{1}{2}$. Since $\pi$ (cuspidal, regular, self-dual) is tempered at all primes, the Archimedean factors $L_\infty(s,\pi)$ and $L_\infty(1-s,\pi)$ are holomorphic at $s=1/2$. Deligne's  conjecture then implies

\vspace{2mm}

\noindent\textbf{Conjecture 1.2.} \textit{For} $a\in \Aut(\C)$
$$
L(1/2,\pi) = 0 \Leftrightarrow L(1/2,a(\pi))=0.
$$

\noindent\textbf{3.} When the field $F$ is totally imaginary, Conjecture 1.2 has been proved by Moeglin \cite{M}. Here we assume $F$ totally real. The Rankin $L$--function $L(s,\pi\times\pi)$ can be written
$$
L(s,\pi\times\pi) = L(s,\pi,\Lambda^2)L(s,\pi,S^2)
$$
associated to the representations $\Lambda^2$ and $S^2$ of the $L$-group. We say that $\pi$ is symplectic (resp. orthogonal) if $L(s,\pi,\Lambda^2)$ $(\resp.~L(s,\pi,S^2)$) has a pole at $s=1$. (This applies for $N$ arbitrary.). For $N$ even and $\pi$ algebraic regular, $\pi$ is always symplectic if $F$ is totally real (Proposition~1.1).

We now impose a stronger condition on $\pi$, superregularity (see Definition~1.1)

Our first result is a new proof of a theorem of Gr\"obner and Raghuram \cite{GR}:

\vspace{2mm}

\noindent\textbf{Theorem 2.5.} \textit{Assume $N$ even and the representation $\pi$ of $\GL(N,\A_F)$  algebraic, self-dual and superregular. Then for any } $a\in \Aut(\C)$,
$$
L(1/2,\pi)=0 \Leftrightarrow L(1/2,a(\pi))=0.
$$

Gr\"obner and Raghuram prove much more, since they obtain formulas for the critical values of $L(s,\pi)$ at all critical  $s\in 1/2+\Z$, from which this follows. However, we will now be able to consider Rankin $L$-functions.

\vspace{2mm}

\noindent\textbf{4.} Assume now $r$ is even, $t$ is odd, and $\pi$ is a superregular, algebraic, self--dual, cuspidal representation of $\GL(r,\A_F)$ while $\rho$ is a regular, algebraic, self--dual, cuspidal of $\GL(t,\A_F)$. We must assume $\pi$, $\rho$ disjoint (Def.~3.1), which is a condition of regularity on their Hodge-types (= representations of $W_\R$ associated to their Archimedean factors.) Recall that $\pi$ is symplectic, while $\rho$ is orthogonal. The conjectural representation $\pi \boxtimes\rho$ (Langlands tensor product) is then symplectic. As in \S~2, we see that the relevant (critical) values will be in $1/2+\Z$. In this case, we prove, for $F$ totally real.

\vspace{2mm}

\noindent\textbf{Theorem 3.2.} \textit{For any} $a\in \Aut(\C)$,
$$
L(1/2,\pi\times\rho)=0 \Leftrightarrow L(1/2,a(\pi) \times a(\rho))=0.
$$
If $F$ is totally imaginary, this is again due to Moeglin.

\vspace{2mm}

\noindent \textbf{Remark.---} Theorem 2.5 follows from Theorem~3.2, by considering a representation $\rho$ of $\GL(1)$.

\vspace{2mm}

\noindent \textbf{5.} The proof relies on Eisenstein cohomology, for the group $\Sp(n)/F$. We summarise the argument for Theorem~2.5. Write $G=\Sp(n)/F$, so $G$ has a Levi subgroup isomorphic to $\GL(n)$. The dual of $G$ is $^L G=\SO(2n+1)\times W_F$. Because $\pi$ is symplectic, there is an Arthur parameter
\begin{equation}
\psi
 = \pi \otimes \usp(2)\oplus \bf{1}
 \end{equation}
 where $\bf{1}$ is the trivial (cuspidal) representation of $\GL(1,\A_F)$. (We refer to Chapter~2 for the details of Arthur's constructions.)
 
 Associated to $\psi$, there is a subspace $\Hh_\psi$ of $L_{disc}^2 (G(F)\ba G(\A_F))$. It can be checked that this space is always non-zero (Corollary to Theorem~2.1). For $\pi$ cohomological, self-dual and superregular, the Archimedean factors of the representations $\pi_G$ occurring in $\Hh_\psi$  are cohomological.
 
 There is a natural candidate for a representation $\pi_G$ in $\Hh_\psi$:  consider, for $s\in\C$, the induced representation
 $$
 I_s=\ind _{MN}^G (\pi \otimes |\det|^s)
 $$
 (induction for the ad\`elic group) where $P=MN$ is the parabolic subgroup associated to $M$. The space of $I_s$ is independent of $s$. For $f\in I_s$, we can form the Eisenstein series $E(f,s)$, defined a priori for $\Re (s) >n-1$ (if $n>1$). Its constant term with respect to $P$ is
 $$
 f\mapsto f+M(s)f
 $$
 where $M(s) : I_s \rg I_{-s}$ has a factor, outside a finite set  $S$ of ramified primes, equal to
$$
\frac{L^S(s,\pi) L^S(2s,\pi,\Lambda^2)}{L^S(s+1,\pi)L^S(2s+1,\pi,\Lambda^2)}.
$$
Since $\pi$ is symplectic, we see that the second factors will contribute a residue at $s=1/2$. It is cancelled if, and only if, $L(1/2,\pi)=0$.

Assume $L(1/2,\pi) \neq 0$. We then obtain a residual representation $\pi_G\subset \Hh_\psi$, with cohomological factors at infinity. Therefore we obtain a summand in the $L^2$-cohomology $\mathbb{H}^q(S_K,\C), S_K$ being the Shimura variety for some level $K$. However, by Zucker's conjecture,
$$
\mathbb{H}^q(S_K,\C) \cong \IH^q(S_K^*,\C)
$$
where the second term is the intersection cohomology of the Baily-Borel compactification. In particular $\mathbb{H}^q(S_K,\C)$ inherits a $\Q$-structure. If $\pi_{G,f}$ is the finite factor of $\pi_G$ and $a\in \Aut(\C/\Q)$, we see that $a(\pi_{G,f})$ occurs in $\mathbb{H}^q(S_K,\C)$:  it comes from a representation $a(\pi_G)$, which is easily seen to belong to $\Hh_{a(\psi)}$ where $a(\psi)$ is defined as in (0.2), $\pi$ being replaced by $a(\pi)$.

Moreover, by a crucial result of Rohlfs and Speh, a form $\omega \in \mathbb{H}^q(S_K,\C)$ --- here $q$ has to be chosen minimal --- restricts non--trivially to the boundary $\partial(S_K^{BS})$ of the Borel-Serre compactification of $S_K$. Since this restriction, seen as $\IH^q(S_K^*,\C) \rg H^q(S_K,\C) \rg H^q(\partial S_K^{BS},\C)$, is defined over $\Q$, we see that, for $a\in \Aut(\C)$, $a(\omega)$ restricts non trivially, and consequently $a(\pi_G)$ is residual.

This, however, does not imply that $a(\pi_G)$ arises, by the formation of residues, from $a(\pi)$. To prove this we must show (replacing $a(\pi)$ by $\pi$) that any residual representation (cohomological)  in $\Hh_\psi$ can only be obtained by the formation of residues of $E(f,s)$ at $s=1/2$, for a function $f$ in the space $I(s)$ associated to $\pi$. This is Theorem~2.3, which relies $(a)$ on a deep result of Grba\v c  and Schwermer, Theorem~2.2, and $(b)$ on a computation using fully Arthur's parametrisation of the discrete representations of $\Sp(m',\A_F)$ in terms of representations of various $\GL(r)$. We refer to the proof of Theorem~2.3.

At this point we know that the residual representation $a(\pi_G)$ can only come from Eisenstein series constructed from $a(\pi)$; the previous computation of the constant term then shows that $L(1/2,a(\pi)) \neq 0$.

We apologise to the reader for the length of our Chapter~2, especially in view of the fact that Theorem~2.5 is not new. However we have found that these arguments were best given first in this simpler case. In Chapter~3, we only have to amplify them to treat Rankin $L$-functions. We now work with the Levi subgroup $M = \GL(r)\times \Sp(n-r)$ of $\Sp(n)$; we are given a cuspidal representation $\pi\otimes \sigma$ of $M(\A_F)$, $\sigma$ being associated (by Arthur's theorems) to a cuspidal, self-dual, orthogonal representation $\rho$ of $\GL(N',\A_F)$ where $N'=2(n-r)+1$. All our data are regular algebraic, so $\rho$ and, therefore, $\sigma$ are tempered. The arguments are then completely similar, once we have proven the requisite extension of Theorem~2.3 --- see Theorem~3.1.

\vspace{2mm}

\noindent \textbf{Caveat.}  Theorems 2.5 and 3.2 both depend on Arthur's book, for which full proofs have not yet appeared. We do not know if proofs could be given using the different approach to functoriality between $\GL(n)$ and classical groups (Kim, Piatestki-Shapiro, Cogdell, Shahidi).

\vspace{2mm}

\noindent \textbf{Acknowledgements.} 
We thank Don Blasius, Harald Grobner, Colette Moeglin, David Renard,  Ga\"{e}tan Chenevier and Claude Sabbah for useful corresondence. As will be clear from the text, we rely heavily on the work of Grbac and Schwermer;  for earlier results  related to ours, we mention again Moeglin and we refer to the extensive list in Raghuram's paper \cite{Ra} \footnote{~page 3}. Finally, we thank the Fondation Simone et Cino Del Duca for material support.

\section{Regular algebraic, self-dual cuspidal representations of $\GL(N)$}

\subsection{}  In this section we consider self--dual cuspidal representations $\pi$ of  $\GL(N,\A)$, where $F$ is a number field and $\A=\A_F$ the ad\`eles of $F$. Let $L(s,\pi)$ be the complete standard $L$-function of~$\pi$:
$$
L(s,\pi) = \prod_v L_v(s,\pi)
$$
where $v$ runs over all places of $F$. We consider also the $L$-functions
$$
L(s,\pi\times \pi) =L(s,\pi,\otimes^2 \St)
$$
(Rankin $L$-function), where $\St$ is the standard representation of $\hat{G}= \GL(N,\C)$; it decomposes as
$$
L(s,\pi\times\pi) = L(s,\pi,S^2) L(s,\pi,\Lambda^2)
$$
where $S^2$ (resp. $\Lambda^2$) denotes the symmetric (antisymmetric) tensors in $\otimes^2 \St$.

Assume now $\pi$ is self-dual: $\pi\cong \tilde{\pi}$. Then $L(s,\pi\times\pi)$ is holomorphic for $\Re(s)>0$,  except for a simple pole at $s=1$. It has been proved by Shahidi and Grba\v c \cite{Gr1} that $L(s,\pi,S^2)$ and $L(s,\pi,\Lambda^2)$ are holomorphic for $\Re(s)>0$, except at most for a simple pole at $s=1$ ; this pole will occur either for the symmetric, or the exterior, square. We say that $\pi$ is \textit{orthogonal} (resp. \textit{symplectic}) as the case may be. If $N$ is odd, $\pi$ is orthogonal: this follows from Arthur \cite[Thm. 1.5.3(a)]{A}, see also Cogdell, Kim, Piatetskii-Shapiro, Shahidi \cite{CKPS1}, \cite{CKPS2}.

If $N$ is even, $\pi$ can be orthogonal or symplectic, as can be seen by considering the case where $\pi$ is associated to an Artin representation $\rho$ of $\Gal(\bar{F}/F)$. The functional equation of $L(s,\pi)$ can be written, $\pi$ being self-dual:
$$
\begin{array}{rl}
L(s,\pi) &=\varepsilon(\pi) D^{1/2-2s}L(1-s,\pi)\\
&=\varepsilon(\pi,s) L(1-s,\pi)
\end{array}
$$
where $D$ is the conductor of $\pi$. This implies that $\varepsilon(\pi,1/2)= \pm 1$.

If $\pi$ is orthogonal, $\varepsilon(\pi,1/2)=1$, cf. \cite[Thm 1.5.3(b)]{A}; earlier results were proved by Lapid \cite{L}. If $\pi$ is symplectic, both signs can generally occur. If $\varepsilon(\pi,1/2)=-1$, $L(1/2,\pi)=0$.

\subsection{} We now consider regular algebraic representations in the sense of \cite{Cl}. Recall the definition. We fix an algebraic closure $\C$ of $\R$. We write $W_\C=\C^\times$, and define $W_\R=\C^\times\coprod j\C^\times$, $j^2=-1, jzj^{-1}=\bar{z}$. (Note that, independently of the choice of $\C$, $W_\R$ is defined up to unique isomorphism.)

Let $v$ be a real prime of $F$. Then $\pi_v$ defines, by the Langlands classification, a representation
$$
r_v : W_\R \lgr \GL(N,\C).
$$

We assume that $r_v |_{\C^\times}$ is of the form
\begin{equation}
z \mapsto ((z/\bar{z})^{p_1},\ldots (z/\bar{z})^{p_N})
\end{equation}
with $p_i\in \Z +\frac{N-1}{2}$, and the $p_i$ are distinct. (Since $\pi$ is self-dual, the weight $w$ of $\pi$ \cite[Lemma 4.9]{Cl} is $0$.) The self-duality further implies that
\begin{equation}
\{p_1,\ldots p_N\} = \{-p_1,\ldots, -p_N\}.
\end{equation}

If now $v$ is a complex prime, we fix an isomorphism $\iota:F_v\cong \C$ and the representation $\pi_v$ of $\GL(N,F_v)$ defines again
$$
r_v : \C^\times \lgr \GL(n,\C)
$$
which must satisfy the same conditions.

In contrast with the case of non-regular representations (e.g., those of Artin type), we now have the following simple result.

\vspace{2mm}

{\maprop Assume $F$ has a real prime $v$. Let $\pi$ be a regular algebraic, self-dual cuspidal representation of $\GL(N,\A)$. Then $\pi$ is orthogonal ($\resp.$ symplectic) if, and only if, $N$ is odd ($\resp.$ even.)}

\vspace{2mm}

This was, in effect, noticed in \cite{CC}, but the formal argument given there can now be made effective. If suffices to consider the case of $N$ even. By Theorem~1.4.1 of \cite{A}, $\pi$ originates from  a unique elliptic simple twisted endoscopic group $G$; then $\hat{G}=\SO(N,\C)$ or $\Sp(N/2,\C) \subset \GL(N,\C)$. By Theorem~1.5.3 of \cite{A}, $\pi$ is symplectic  (in our sense) if, and only if, $\hat{G}=\Sp(N/2,\C)$.

If $\hat{G}$ is orthogonal, $G/F$ is isomorphic to the split group $\SO (N)$ or to the quasi-split, non-split $\SO^*(N)$ split over a quadratic extension $F'/F$. In the first case there is an obvious morphism $^LG \lgr \GL(N,\C)$. In the second one, we get a priori $(^LG)_{F'/F}\lgr (^L\GL(N))_{F'/F}$ (see Arthur \cite[\ p. 10--11]{A}) but this can be seen as a homomorphism $\uO(N,\C)\lgr \GL(N,\C)$ by the choice of a suitable element in $\uO(N,\C)$ of determinant $(-1)$.

Assume $\pi$ is orthogonal and consider the Langlands parameter $r_v : W_{F_v}\rg \GL(N,\C)$ associated to $\pi_v$. This quotients \cite[Thm.~1.4.2]{A} through a map $W_{F_v} \lgr \uO(N,\C) \lgr \GL(N,\C)$. In  particular we get an orthogonal representation of $W_{F_v}=W_\R$. However, by (1.1) and (1.2) the representation associated to $\pi_v$ is the following. Let $\chi_i(z) =(z/\bar{z})^{p_i}$ and $r_i=\ind_{\C^\times}^{W_\R} \chi_i$ such that $r_v = \bigoplus r_i$. The representation $r_i$ is irreducible and self-dual. Its determinant, a character of $W_\R$, thus of $\R^{\times}$, is by a well-known formula for ``transfer''
$$
\det r_i = \chi_i|_{\R^\times} \varepsilon_{\C/\R}
$$
where $\varepsilon_{\C/\R}$ is the sign character of $\R^\times$. Since $p_i \in 1/2 +\Z$, $\det r_i=1$ and the image of $r_i$ is contained in $\SL(2,\C):r_i$ is symplectic, and the $\chi_i^{\pm 1}$ being distinct any invariant form on $\C^N$ is symplectic. This contradicts the fact that $r_v$ is orthogonal.

\subsection{} Recall now the arithmeticity  properties of regular algebraic representations, cf. \cite{Cl}. We assume again that $F$ is a general number field. For regular algebraic representations --- not necessarily self-dual --- condition (1.1) is replaced by
\begin{equation}
z \longmapsto (z^{p_1}\, \bar{z}^{q_1}, \ldots z^{p_N} \bar{z}^{q_N})
\end{equation}
with $p_i+q_i=w$, and $p_i \in \Z +\frac{N-1}{2}$.

Let $\pi_f = \bigotimes\limits_{v\nmid \infty} \pi_v$. It is then known that $\pi_f$ is defined over a number field $E\subset \C$. We can take $E$ to be Galois over $\Q$ ; then, for any $a\in \Gal(E/\Q)$, $a(\pi_f)=\pi_f \otimes_{E,a}E$ is the finite part of a cuspidal, regular algebraic representation $a(\pi)$. (The type of $a(\pi)$ at infinity can be described, cf. \cite[Th\'eor\`eme~3.13]{Cl}\footnote{The statement of Thm. 3.13 is incorrect: $\pi$ should be assumed cuspidal.}.)

(In these arguments, one may also consider $a(\pi_f)$ for $a\in \Aut(\C/\Q)$, with the same arguments, thus defining $a(\pi)$.) 

\subsection{} 
Now fix $F$, and consider absolutely irreducible motives $M$ over $F$, with coefficients in a number field $E$ \cite[\S 1-2]{Del}.

Given an embedding $\sigma: E \rightarrow \C$, one can form the $L$-function $L(\sigma,M,s)$ \cite[\S2.2]{Del}. It can be completed by an Archimedean factor $L_{\infty}(M,s)$, independent of $\sigma$. Let
$$
\Lambda(\sigma,M,s)=L_{\infty}(M,s) L( \sigma,M,s).
$$
It is expected that $\Lambda(\sigma,M,s)$ extends meromorphically to $\C$, with functional equation
$$
\Lambda(\sigma,M,s)= \varepsilon(\sigma, M,s) \Lambda(\sigma, \tilde{M},s)  
$$
where $\tilde{M}$ is the dual motive. An integer $n$ is critical for $M$ if $L_{\infty}(M,s)$ and  $L_{\infty}(\tilde{M},s)$  are holomorphic at $n$. Deligne then conjectures that $L(\sigma,M,n)$ is given by a period of the motive $M$. This implies:

\vskip3mm

\textit{For} $\sigma, \sigma' :E \rightarrow \C, L(\sigma,M,n)=0  \Leftrightarrow L(\sigma',M,n)=0$ \textit{if} $n$ \textit{is critical}.





 \vskip3mm

 It is also conjectured that $L(\sigma,M,s)= L(s,\pi)$ for an irreducible cuspidal (complex) representation $\pi$ of $\GL(N,\A_F)$, where $N$ is the dimension of $M$. At an Archimedean place $v$ of $F$, the representation of $\C^{\times}$ associated to $\pi_v$ is deduced from the Hodge structure on $M_{Betti} \otimes \C$, where the tensor product is given by $\iota: E \rightarrow \C$ defining $v$. In particular, if $N$ is odd $\pi$ is algebraic. If $N$ is even, this is not the case.
 
 \vskip3mm
 
 Assume $N$ is even, and consider $\pi'=\pi[-1/2] := \pi \otimes \vert \rm det\vert ^{-1/2}$ where $\vert ~\vert$ is the id\`{e}le norm. It is algebraic and $\tilde{\pi}' = \tilde{\pi} [1/2]$. One checks easily that
 
 \vskip3mm
 $n$ \textit{ critical for} $\pi \Leftrightarrow L_{\infty}(\nu, \pi')$ \textit{and}  $L_{\infty}(1 - \nu, \pi')$ \textit{holomorphic}, $\nu=n+ 1/2$. 
 
 \vskip3mm
 
 In this case, we say that $\nu$ is critical for $\pi'$. On the other hand, if $a\in \Aut(\C)$,
 $$
 L(\sigma,M,s)=L(s,\pi)= \sum c_n n^{-s} 
 $$
 with $c_n \in \sigma(E) \subset \C$ . Thus
 $$
 L(s+1/2, \pi')= \sum c_n n^{-s}
 $$
 and
 $$
 L(s+1/2,a( \pi')) = \sum a(c_n) n^{-s}
 $$
 by the computation in \cite[\S 4.3.4]{Cl}. \footnote{ In \cite{Cl} the translation of $s$ is by $\frac{N-1}{2}$ rather than $1/2$, but this does not change the argument.} This is $L(a\sigma,M,s)$. Since we expect conversely that $\pi'$ will also originate in a motive $M'$, we see that
 
 {\maconj For $\nu \in 1/2+\Z$ critical,  $\pi$ algebraic, 
 $$
 L(\nu, \pi) = 0 \Leftrightarrow  L(\nu, a(\pi)) = 0.
 $$}
 
 
 Assume in particular $N$ even and $\pi$ algebraic, cuspidal, regular and self-dual. Then $1/2$ is critical since $\pi_{\infty}$ is tempered and its $L$-function is holomorphic for $\Re(s) \in ]0, 1]$. Thus
 
 {\maconj For $a \in \Aut(\C/\Q)$,
  $$
  L(1/2, \pi)=0 \Leftrightarrow L(1/2, a(\pi))=0.
  $$}

\section{The case of symplectic representations}

\subsection{} 
In this chapter we assume ${n}=2m$ even and consider a regular algebraic, self-dual, cuspidal $(RASDC$ for short) representation $\pi$ of $\GL(n,\A_F)$. We assume $F$ totally real. We recall that the Ramanujan conjecture is then known for $\pi$: all factors $\pi_v$ are tempered. Cf. \cite{Ca}, \cite{ CP}.\footnote{For all the arguments in this paper, the easier result of \cite{CP} suffices.} Moreover $\pi$ is symplectic.

Let $G=\Sp(n) /F$. Then $\hat{G}=\SO(2n+1,\C)$. We consider an Arthur parameter
$$
\psi=\pi \otimes \usp(2) \oplus \bf{1}
$$
(Arthur \cite[1.4]{A} where $\usp(2)$ is the $2$-dimensional representation of $\SL(2,\C)$. Here $\bf{1}$ is the trivial (cuspidal) representation of $\GL(1,\A_f)$. The group $\Ll_\psi$ associated to $\psi$ (Arthur, loc. cit.) is $\Sp(m)\times \SO(1)=\Sp(m)$. There is a natural embedding
$$
\begin{array}{l}
\Ll_\psi \times \SL(2,\C) \lgr \hat{G}\\ \vspace{2mm}
\quad (\ell,u) \longmapsto \ell \otimes \usp_2(u) \oplus 1.
\end{array}
$$
Since $\ell$ is symplectic $\ell \otimes \usp_2(u)$ is orthogonal.

This can be extended to a homomorphism of split $L$-groups. Since $\pi$ is symplectic, $\pi_v$ defines a Hecke matrix $t_v\in \Sp(n/2,\C)$ for all unramified primes. Then
$$
t_v \otimes \Big(
\begin{array}{ll}
q_v^{1/2} \\
&\ q_v^{-1/2}
\end{array}
\Big) \oplus 1
$$
is a Hecke matrix (in $\SO(2n+1,\C))$ for $G$. The parameter $\psi$ belongs to $\tilde{\Psi}_2(G)$ (Arthur \cite[p.~33]{A}); Theorem~1.5.2 of \cite{A} states that it contributes to $L_{disc}^2 (G(F)\ba G(\A)$. Precisely, for each place $v$, $\psi$ defines a local parameter $\psi_v:L_{F_v}\times \SL(2,\C) \rg \hat{G}$ where 
$$
\begin{array}{rl}
L_{F_v} &= W_{F_v} \quad (v\ \mathrm{Archimedean})\\
L_{F_v} &= W_{F_v}\times \SL(2,\C) \quad (v\ \mathrm{finite})
\end{array}
$$
Cf. \cite[p. 40]{A}.

There is a group $\CSs_\psi$ associated to $\psi$, the centraliser of $\Ll_\psi\times \SL(2,\C)$  in $\hat{G}$. It is equal to $(\pm1)$. For each $v$, there is a finite set $\tilde{\Pi}_{\psi_v}$ of representations-with-multiplicities of $G(F_v)$ and a pairing
$$
\pi_G \longmapsto \langle \cdot , \pi_G \rangle \quad (\cdot \in \CSs_\psi, \pi_G \in \tilde{\Pi}_{\psi_v}).
$$
\cite[Thm. 1.5.1]{A}.

Let $\Hh_\psi \subset L^2_{disc}(G(F)\ba G(\A))$ be the subspace associated to $\psi$. Then
$$
\Hh_\psi \cong \bigoplus_{\pi_G \in \tilde{\Pi}_\psi(\varepsilon_\psi)}\pi_{G}
$$
where $\tilde{\Pi}_\psi$ is the (restricted) tensor product of the representations in $\tilde{\Pi}_{\psi_v}$ --- see \cite[p.~47]{A} for details --- and $\tilde{\pi}_\psi(\varepsilon_\psi)$ is defined by
$$
\prod_v \langle s,\pi_{G,v} \rangle = \varepsilon_\psi(s) \quad (s\in \CSs_\psi)
$$
where $\varepsilon_\psi : \CSs_\psi \rg \pm 1$ is a character.

To understand $\Hh_\psi$, we must first compute~$\varepsilon_\psi$.

\vspace{2mm}

{\maprop Let $s_0\in \CSs_\psi$ be the non-trivial element. Then}
$$
\varepsilon_\psi(s_0) = \varepsilon(1/2,\pi)\ (=\pm 1).
$$

The character $\varepsilon_\psi$ is defined by the adjoint action $\tau_\psi$ of $\CSs_\psi \times \Ll_\psi \times \SL(2,\C)$ on $\hat{\gog}=\Lie\, \hat{G}$,~so
\begin{equation}
\tau_\psi = \bigoplus_\alpha \lambda_\alpha \otimes \mu_\alpha\otimes \nu_\alpha,
\end{equation}
decomposition into irreducibles for $\CSs_\psi$, $\Ll_\psi$ and $\SL(2,\C)$. Here $\Ll_\psi = \Sp(n/2,\C) $ embeds  in a Levi subgroup $\hat{M}=\GL(n,\C)$ in $\hat{G}$. We have the  triangular decomposition $\hat{\gog}= \hat{\fm}\oplus \hat{\fn}\oplus \hat{\fu}$; as a representation of $\hat{M}$, $\hat{\fn} = \St_n \oplus \Lambda^2$, $\hat{\fu}=(\St_n)^* \oplus (\Lambda^2)^*$. The representation $\St_n$ of $\Sp(n/2,\C)$ is self-dual and so is~$\Lambda^2$.

In the decomposition (2.1) the irreducibles $\mu_\alpha$, representations of $\Sp(n/2,\C)$, are self-dual; we are only interested in those that are symplectic; $\mu_\alpha$ (with the data $L_{F_v}\rg \Ll_\psi)$ defines a global $L$-function $L(s,\mu_\alpha)$ and an $\varepsilon$-factor $\varepsilon(1/2,\mu_\alpha)$. See \cite[p.~48]{A}.

Recall that the irreducible representations of $\Sp(g,\C)$ (with its usual diagonal torus $\cong(\C^\times)^g)$ are parametrized by their highest weight $\xi \in \Z^g$, with $\xi_1\ge \xi_2\ge \cdots \xi_g \ge 0$. The following lemma is probably well-known, but we have not found a reference.

\vspace{2mm}

{\monlem Let $V=V^\xi$ be the irreducible representation of $\Sp(g,\C)$ with highest weight $\xi$. Then $V$ is orthogonal ($\resp.$ symplectic) iff $\sum \ \xi_i\equiv 0$ $[2]$, $\resp. ~\sum\ \xi_i \equiv 1$ $[2]$. }

\vspace{2mm}

The group $\Sp(g)$ naturally contains $\Sp(1)^g=H$ (associated to the long roots). Assume for instance that $V$ is symplectic. The irreducible representation $W$ of $H$ of highest weight $\xi$ is a submodule of $V$. It occurs with multiplicity 1 since $\xi$ (as a character of the maximal torus $T_H=T=(\C^*)^g)$ does. Thus the invariant bilinear form $B$ on $V$ is non-zero on $W$, since $W$ is self-dual and occurs with multiplicity one. 

We can write $W=W_1 \otimes \cdots \otimes W_g$, $W_i$ having highest weight $\xi_i$. Then $W_i$ is orthogonal/symplectic according to whether $\xi_i \equiv 0/\xi_i\equiv 1$, and $W$ is symplectic if, and only if, $\sum \,\xi_i\equiv 1$.

\vspace{2mm}

The lemma implies that, in the foregoing decompositions, only $\St_n\subset \hat{\fn}$ and $\St_n^* \cong \St_n\subset \hat{\fu}$ contribute symplectic representations. Looking at the triangular decomposition of $\hat{\gog}$, one checks that the action of $\Ll_\psi \times \SL(2,\C)$ is isomorphic to $\St_n\otimes \usp(2)$. On the other hand, we have $s=\Big(\begin{smallmatrix}
s_1\\
&1\\
&&s_1
\end{smallmatrix}\Big)$ compatibly with the triangular decomposition, with $s_1=\pm1$; then the adjoint action of $s$ on the summands $\St_n$ and $\St_n^*$ is by $X \mapsto s_1X$. We deduce that (the only relevant) $\lambda\otimes \mu \otimes \nu$ is $\sgn \otimes \St_n \otimes \usp(2)$; Arthur's definition of $\varepsilon_\psi$ is then:
$$
\begin{array}{rlcl}
\varepsilon
_\psi(s) &=\sgn &\mathrm{if} &\varepsilon(1/2,\pi)=-1\\
&\equiv 1 &\mathrm{if} &\varepsilon(1/2,\pi)=1,
\end{array}
$$
recalling that the factor $\St_n$ defines the $L$-function associated to $\mu=\pi$.

\subsection{} The description of the data $\tilde{\Pi}_{\psi_v}$ associated to $\psi$ is in general (for our non-tempered $\psi$) difficult. It can be made more explicit, at the Archimedean primes, in the cohomological case.

A still indirect, but more concrete, description of $\tilde{\Pi}_{\psi_v}$ is given by Arthur in \cite{A2}. The global ``packet'' $\tilde{\Pi}_{\psi}$ is related to the following global representation of $\GL(N,\A)$, $N=2n+1$. We write $\ind_{a_1,a_2}^a$ $(a=a_1+a_2)$ for the unitary induction to $\GL(a)$ of a representation of $\GL(a_1)\times \GL(a_2)$. We write $\pi[s]$ for $\pi\otimes |\det|^s$. First, $\ind_{n,n}^{2n}(\pi[1/2]\otimes \pi[-1/2])$  has a unitary quotient $\Ll(\pi[1/2],\pi[-1/2])$ (Moeglin-Waldspurger). Then we consider the ad\`elic representation
$$
\Pi= \ind_{2n,1}^{2n+1}(\Ll(\pi[1/2],\pi[-1/2])\otimes 1).
$$
This also defines $\Pi_v$ at any prime; it is a unitary representation.

The representation $\Pi_v$ is $\theta$-stable, for the automorphism of $\GL(N)$ realising the  duals of representations; $\theta$ has to be chosen suitably, cf. Arthur \cite[p.~7--9]{A}; in particular $\theta$ leaves the Borel subgroup invariant. There is then a canonical intertwining operator $I_\theta:\Pi_v \rg \Pi_v \circ \theta$. (It is defined by the theory of Whittaker models.)

For a function $\varphi$ on $\GL(N,F_v)$ we consider $\tilde{\varphi}(\psi) = \trace (\Pi_v(\varphi\,)I_\theta)$.

One can define associated functions $(f,\varphi)$ where $f$ is (smooth compactly supported) on $G(F_v)$ and $\varphi$ on $\GL(N,F_v)$. If $\rho$ is a representation of $G(F_v)$ let $f(\rho) = {\rm trace}~ \rho(f)$. Then, for $(\varphi,f)$ associated:
$$
\tilde{\varphi}(\psi) = \sum_{\rho\in\tilde{\Pi}_\psi} \langle s_\psi,\rho\rangle f(\rho) 
$$
where $s_\psi$ is the image of  $\Big(1,(\begin{smallmatrix}
-1\\
&-1
\end{smallmatrix})\Big) \in L_F \times \SL(2,\C)$ by the local version of $\psi:L_{F_v} \times \SL(2,\C) \rg \GL(N,\C)$. For details we refer to \cite[\S 30]{A2}.

Since $f$ and $\varphi$ are associated by their (twisted) orbital integrals, this implies at an  Archimedean prime that the infinitesimal  characters of $\Pi_v$ and of all representations $\rho \in\tilde{\Pi}_{\psi_v}$ are associated. We proceed to describe this correspondence.

Choose a split maximal torus $T\subset G(F_v)$. Then $\Lie\, T\otimes \C\ \cong \C^n , \Lie\,\hat{T} \otimes \C \cong \C^n$.

Similarly if $T_N\subset \GL(N)$ is a maximal split torus, $\Lie\, T_N \otimes \C \cong \C^N, \Lie\,\hat{T_N} \otimes \C \cong \C^N$.

An infinitesimal character for $G(\R)$ is  then $X=(X_1,\ldots X_n)$; an infinitesimal character of $\GL(N,\R)$ is
$$
Y=(Y_1,\ldots Y_N).
$$
The representation $\Pi_v$ is self-dual, and this implies $Y=-Y$ $(\!\!\!\mod \,\Ss_N)$, so we can write
$$
Y=(Y_1,\ldots,Y_n,0,-Y_n,\ldots -Y_1).
$$
Then $X=(Y_1,\ldots,Y_n)$; $X$ is well-defined modulo $W_G=(\pm1)^n \rtimes \Ss_n$. Given the Langlands parameters of $\pi$ (\S~1.1) we see that $(\!\!\!\mod W_G)$  $X$ is equal to
$$
(p_1 +1/2,\ p_1-1/2,\ p_2+1/2,\ldots, p_m+1/2,\ p_m-1/2)
$$
with $p_i\in \frac{1}{2}+\Z$, $p_1>\cdots> p_m > 0$.

In particular $p_i\pm 1/2$ is now an integer, and $X$ is a weight of $G$ in the dual of $\Lie(T)$. It is strictly dominant if
$$
\begin{array}{lrcl}
&p_1+1/2 &>&p_1-1/2 > \cdots\ \cdots > p_m-\frac{1}{2},\\
i.e. &p_i &\ge &p_{i+1} +2 \qquad (i=1,\ldots m-1)\\
&p_m &\ge &\frac{3}{2}.
\end{array}
$$

\vspace{2mm}

{\madef Assume $\pi$ is a regular algebraic, self dual representation of $\GL(n,\R)$, $n=2m$. Let $p=(p_1,\ldots p_n)$ be its infinitesimal character,
$$
p=(p_1,\ldots p_m,-p_m,\ldots -p_1)\quad p_1>p_2>\cdots >-p_1.
$$
We say that $\pi$ is \rm {superregular} \textit{if}, $\mod \Ss_n$,}
$$
\begin{array}{rcl}
p_i & \ge& p_{i+1}+2 \qquad (i=1,\ldots m-1)\\
p_m&\ge &\frac{3}{2}.
\end{array}
$$

If $\pi_G$ belongs to the $\psi$-packet associated to $\psi_v$, we see that the infinitesimal character of $\pi_G$ is regular integral. Since $\pi_G$ is unitary, this implies that $\pi_G$ is a cohomological representation of $G(F_v)=G(\R)$ \cite{Sa}. We record the result:

\vspace{2mm}

{\maprop Assume $\pi_v$ is superregular at all real primes. Then, for each such prime $v$, all representations in $\tilde{\Pi}_{\psi_v}$ are cohomological (unitary).}

\vspace{2mm}

Consider, for given $v$, the sign $s\mapsto \langle s,\pi_G\rangle$ for $\pi_G\in \tilde{\Pi}_{\psi_v}$. Moeglin and Renard have studied these Arthur packets \cite{MR}, and their results imply inter alia:

\vspace{2mm}

{\montheo For $\psi$ as above, $v$ real and $\pi$ $(RASDC)$ superregular, the sign $s\mapsto \langle s,\pi_G \rangle$ is the trivial ($\resp.$ non-trivial) character of $\CSs_\psi$ for some $\pi_G \in \tilde{\Pi}_{\psi_v}$.}

\vspace{2mm}

\noindent\textbf{Corollary.} \textit{The space $\Hh_\psi \subset L_{disc}^2(G(F)\ba G(\A))$ associated to $\psi$ is non-zero.}

\vspace{2mm}

This is clear: we must fulfill the condition
$$
\prod_v \langle s_0,\pi_{G,v} \rangle = \varepsilon_\psi(s_0)
$$
where $s_0$ is the generator of $\CSs_\psi$. Assume for instance $\varepsilon_\psi(s_0) =1$ (cf. Proposition~2.1). The product over finite primes takes the value $1$ or $-1$ for a choice of the $(\pi_{G,v})_{v\nmid \infty}$. We then choose the $\pi_v$ $(v|\infty)$ so as to make the total product equal~to~$1$.

\vspace{2mm}
[We sketch the proof of the theorem, which was explained to us by David Renard. At a real prime, the parameter $\psi_v$ is equal to 
$$
\psi_v=\bigoplus_{i=1}^m V(p_i,-p_i) \otimes \usp(2) \oplus \bf{1}
$$
where $V(p_i,-p_i)$ is the discrete series representation of $\GL(2,\R)$ associated to the representation $\ind_{\C^{\times}}^{W_{\R}}(z/\bar{z})^{p_i}$ . The corresponding (local) $\CSs$-group is $\CSs_{\psi_v}= (\Z/2\Z)^m$;   $\CSs_{\psi_v}$ contains $\CSs_{\psi}  $ embedded diagonally. We can write
$$
\psi_v= V(p_1,-p_1) \oplus \psi'
$$
$\psi'$ being a parameter for $\Sp(n-1,\R) := G'$; obviously    $\CSs_{\psi'} =( \Z/2\Z)^{n-1}$; we assume that the Arthur packet $\tilde{\Pi}_{\psi'}$ is defined; it is composed of cohomological (unitary) representations. The construction of Moeglin-Renard \cite {MR} realises $\tilde{\Pi}_{\psi}$ by cohomological induction from a $\theta$-stable  parabolic algebra $\frak{q}$ whose Levi component $\frak{l} \subset \frak{g}$ is associated to a complex subgroup $L(\C)= \GL(2,\C) \times \Sp(n-1,\C) \subset \Sp(n,\C)$. One considers groups $L_i$ that are pure inner forms \cite{MR} of $ \uU(2) \times \Sp(n-1, \R)$; i.e., 
$$
L_i= \uU(2, 2-i) \times \Sp(n-1,\R), i=0,1,2.  
$$

Let $\xi_i$ be the (Abelian) character of $U(i,2-i)$ of differential $(p_1-1/2, p_1-1/2)$  \cite [p.1861]{MR}. Then $\tilde{\Pi}_{\psi}$  is composed of the representations of $G$ obtained by cohomological induction:
$$
\pi_G= R_{\frak{q}}^d(\xi_i \otimes \pi')
$$ 
and $d$ is defined in \cite[(5.4.4)]{MR}.

The sign $\langle s_0, \pi_G\rangle$ associated to $\pi_G$ is the product of the sign associated to $\pi'$ with $(-1)^i$, cf. \cite[7.1.2]{MR}. (The constant $c$ in this formula is defined p. 1830; here c=1.) For varying $i$ we get all signs.]

\subsection{} With the parameter $\psi$ fixed as above, consider an irreducible global representation $\pi_G \subset \Hh_\psi$. It is automorphic, and therefore obtained from \mbox{(multi-)re\-si\-dues}  of  Eisenstein series formed from parabolic subgroups of~$G$.

At almost all primes $v$, the representation $\pi_G$ is unramified and its Hecke matrix (in $\hat{G}=\SO(2n+1,\C)$ ) is defined by the parameter $\psi$. We have $\hat{T}=\lbrace\mathrm{diag}(t_1,\ldots t_n,1,t_n^{-1},\ldots t_1^{-1})\rbrace \subset \GL(N,\C)$. The element $t=(t_1,\ldots t_n)$ is then uniquely defined $\mod W_G=(\pm 1)^n \rtimes \Ss_n$.

At an unramified prime $v$, the Hecke matrix can be taken to be
\begin{equation}
t_{\pi_G,v} = t_{\psi,v} = t_\pi q^{1/2}
\end{equation}
where $q=q_v$.

The representation $\pi_G$ can be cuspidal or can be obtained from (residues and principal values of) Eisenstein series induced from a non-trivial parabolic subgroup. A standard Levi subgroup of $G$ is of the form $M=\GL(n_1)\times\dots \times \GL(n_r)\times G'$, where $G' = \Sp(n-r)$. 

We consider a datum $(\pi_1,\ldots\pi_r,\sigma)$ where $\pi_i$ ($\resp. \sigma)$ is a cuspidal representation of $\GL(n_i,\A)$ $(\resp. G'(\A))$. We consider the induced representation
$$
\ind_{MN(\A)}^{G(\A)} (\pi_1[s_1] \otimes \cdots \otimes \pi_r [s_r] \otimes \sigma)
$$
which  we write $\pi_1[s_1] \oplus \dots  \oplus\pi_r[s_r] \oplus \sigma$. Its Hecke matrix at $v$ is then
$$
t_1 q^{s_1} \oplus \cdots \oplus t_r q^{s_r} \oplus t_1^{-1} q^{-s_1} \oplus \cdots \oplus t_r^{-1} q^{-s_r} \oplus t_\sigma
$$
where $t_i = t_{\pi_i}$ and $t_{\sigma}\in \hat{G}'$ is seen as before as an element of $T_{\hat{G}'}$. It is defined modulo $W_G$. We write $\pi_M=(\pi_i,\sigma)$. (If $ G'$ is trivial, $ t_{\sigma}=1$ and the following arguments also apply.)

\vspace{2mm}

{\montheo {(Grba\v c, Schwermer)}. Assume a representation $\pi_G$, obtained from multi-residues and principal values of Eisenstein series from $(M,\pi_M)$ at the point $s = (s_1, \ldots s_r)$, is cohomological. Assume $\pi_i|_{\R_+^\times}=1$ where $\R_+^\times$ is embedded diagonally in $\GL(n_i,F_{\infty})$. Then $s_i\in 1/2\  \Z$.}

\vspace{2mm}

See \cite[Theorem 4.1]{GS1} as well as the exposition in \cite[Theorem 3.2]{Gr2}. (We note that the statement in \cite{Gr2} is imprecise, as it assumes only the $\pi_i$ to be unitary:  replacing $\pi_i$ by $\pi_i[it]$ would replace $s$ by $s-it$. The correct condition on $\pi_i$ is stated in \cite[Proposition 1.1]{GS1}.)

We note that there is an obvious candidate for a representation $\pi_G\subset \Hh_\psi$. Namely, consider $M=\GL(n)\subset G$, and the induced representation
$$
\pi_{G,s} = \ind_{MN}^G(\pi[s])
$$ 
(ad\`elic  induced representation.) Its Hecke matrix is given by (2.2) if $s=1/2$.  If it gives rise to a residual representation, this representation will belong to $\Hh_\psi$. We will return to this in the next paragraph.

\vspace{2mm}

{\montheo Assume $\pi_G\subset \Hh_\psi$ is not cuspidal, is cohomological and obtained from multi-residues and principal values of Eisenstein series from $(M,\pi_M)$ at a point $s$. Then $M=\GL(n)$, $\pi_M=\pi$ (up to association) and $s=1/2$.}

\vspace{2mm}

By Arthur's results, $\sigma$ is associated to an Arthur parameter of dimension $2n-2r+1$ 
$$
\psi' = \bigoplus \sigma_j \otimes \usp(m_j)
$$
where $\sigma_j$ is a self-dual, cuspidal representation of $\GL(\nu_j)$ and $\sum \nu_j m_j=2n-2r+1$.

Fix a finite place $v$ where all data are unramified. The corresponding Hecke matrix, seen as a matrix of  dimension  $2(n-r)+1,$ is
$$
\bigoplus t_{\sigma_j} \otimes \usp(m_j)
$$
where we also write $\usp(m_j)$ for $\mathrm{diag}(q^{\frac{m_j-1}{2}},\ldots,q^{\frac{1-m_j}{2}})$. The Hecke matrix of $\ind _{MN}^G(\pi_1[s_1]\otimes\cdots\otimes \pi_r[s_r]\otimes \sigma)$ is then
\begin{equation}
T=\bigoplus_i (t_{\pi_i} q^{s_i} \oplus t_{\tilde{\pi}_i}q^{-s_i}) \oplus \bigoplus_j (t_{\sigma_j} \otimes \usp(m_j)),
\end{equation}
written again as a matrix of dimension $2n+1$. We compare this with (2.2), in dimension~$2n+1$:
\begin{equation}
T'=t_\pi \, q^{1/2} \oplus t_\pi q^{-1/2} \oplus 1.
\end{equation}
(Recall that $\pi$ is self-dual.) Now $\pi$ is cohomological, self dual. For any of the representations $\pi$, $\pi_i$, $\sigma_j$ we write $t^\alpha$ for the diagonal entries of the Hecke matrix. By purity \cite{CP} we have, for $\pi$:
$$
|t_\pi^{\alpha}|=1.
$$

The representations $\pi_i$ are \textbf{unitary} cuspidal, as are the $\sigma_j$. The entries $t^\alpha$, in each case, satisfy
\begin{equation}
q^{-1/2+\delta} <|t^\alpha| < q^{1/2-\delta}
\end{equation}
where $\delta, 0<\delta <1/2$ ,is given by the theorem of Luo-Rudnick-Sarnak \cite{LRS}. Comparing (2.3) and (2.4), we see in the case of $\sigma_j$ that
$$
q^{1/2} = |t^\alpha| q^{\frac{m_j-1}{2}}
$$
whence by (2.5) $m_j\le 2$, or that
$$
1=|t^\alpha| q^{\frac{m_j-1}{2}}, \ \mathrm{whence\ }m_j=1.
$$
Moreover $(m_j=1)$ can occur only once. It does occur since $2(n-r)+1$ is odd.

We now consider the eigenvalues $t^\alpha q^{s_i}$, $t^{-\alpha }q^{-s_i}$, $t^\alpha$ eigenvalue of $\pi_i$, with $s_i\in \frac{1}{2}\Z$. Let $\lambda= \log_q |t^\alpha|$. If $|t^\alpha q^{s_i}|=q^{1/2}$, $\lambda+s_i=1/2$ with $|\lambda|< 1/2-\delta$. This implies $s_i=1/2$ and $\lambda=0$. Similarly, $t^\alpha q^{s_i}=q^{-1/2}$ implies $\lambda=0$, $s_i=-1/2$.
(We do not have to consider the eigenvalue $1$ of (2.4), as it is accounted for by $\sigma$.)

Finally, we see that $s_i=\pm1/2$ for all $i$; however by the invariance of the induced representation by association, we can assume $s_i=1/2$. We also note that the argument implies that the $\pi_i$ are tempered (at all unramified primes). Choosing representatives in $\hat{T}$, we see that
$$
\bigoplus_i (t_{\pi_i }q^{1/2}) \oplus \bigoplus _j(t_{ \sigma_j }q^{1/2}),
$$
the second sum ranging over such $j$'s that $m_j=2$, is equal (i.e., conjugate in $\Ss_n$) to $t_\pi q^{1/2}$. This implies that the representation
$$
\pi_1 \boxplus\cdots \boxplus \pi_r \boxplus \sigma_1 \boxplus \cdots \boxplus \sigma_s
$$
(induced representation of $\GL(n,\A)$) has the same Hecke matrix as $\pi$, at almost all primes $v$. (Moreover $\sigma_j = \bf{1}$ for the unique $j$ such that $m_j=1$). By the results of Jacquet-Shalika \cite{JS}, this implies that $\pi=\pi_1$ or $\pi=\sigma_1$ (after reindexing). In the second case, $\sigma$ is cuspidal by assumption, $G'=G$ and $\pi_G$ is cuspidal. In the first case, $\pi = \pi_1, M = \GL(n)$ and we obtain the assertion of the Theorem.

\vspace{2mm}

\noindent \textbf{Note:} This is a strengthening of Theorem~3.1 of Grba\v c \cite{Gr1}; Grba\v c's theorem describes only the part of $\Hh_\psi$ obtained from Eisenstein series from~$P=MN$.

\subsection{} We now assume that there is a residual representation in $\Hh_\psi$; it is then obtained from Eisenstein series $E(f,s)$ where $f$ belongs to the space, independent of $s$,
$$
I=\ind _{P=MN}^G (\pi[s])=I(s)
$$
with $M=\GL(n)$. The corresponding constant term, for $f\in I$ unramified outside $S\supset \infty$, has a factor (outside $S$) equal to
$$
\frac{L^S(s,\pi) L^S(2s,\pi,\Lambda^2)}{L^S(s+1,\pi)L^S(2s+1,\pi,\Lambda^2)}.
$$
In particular we obtain residues at $s=1/2$ if, and only~if
\begin{equation}
L(\frac{1}{2},\pi) \neq 0
\end{equation}
(as the factors for $v\in S$ are convergent and $\not=0$). See Grba\v c \cite{Gr1}. In this case, $\Res_{s=1/2}$ yields an intertwining  operator
$$
R(1/2) : I(\frac{1}{2}) \lgr L_{disc}^2 (G(F)\ba G(\A).
$$
It image has been described by Grba\v c and Schwermer\footnote{\cite[ p. 14--15, and Thm 8.2]{GS2}}. Note that $I_v(s)$ is for any $s$ with $\Re(s)>0$ an induced representation with $\pi[s]$ positive in Langlands's sense. It has a unique irreducible Langlands quotient $J_v(s)$ and the image of $R(1/2)$ is
$$
\bigotimes_v J_v(1/2).
$$

We write $\Ll(\pi_v,1/2)$ for this Langlands quotient. For $v|\infty$, $\Ll(\pi_v,1/2)$ is unitary, with regular integral infinitesimal character, and therefore cohomological. For simplicity of redaction we first assume that it has cohomology  with trivial coefficients. (The additions in the general case will be explained at the end of this chapter.) Let $q_v$ be the smallest degree $i$ such that
$$
H^i(\gog_v,K_v ; \Ll(\pi_v,1/2)) \neq 0
$$
with obvious notation. (The value of $q_v$ is irrelevant to us; it has been computed by Grba\v c  and Schwermer: see \cite [Theorem 8.3]{GS2}. Then $\dim H^{q_v}(\gog_v, K_v,\break \Ll(\pi_v,1/2))=1$, cf.~\cite{VZ}.

We write, for $K\subset G(\A_f)$ compact open, $S_K=G(F)\ba G(\A)/K_\infty K$, $G_K=G(F)\ba G(\A)/K$, with $K_\infty \subset G(F_\infty)$ maximal compact. Consider the cohomology of $S_K$. We have first the cohomology of the discrete spectrum:
\begin{equation}
\mathbb{H}^\bullet(S_K) = H^\bullet (\gog,K_\infty; L_{disc}^2 (G_K));
\end{equation}
this is the space of $L^2$-harmonic forms on $S_K$, equal to the $L^2$-cohomology of $S_K$ in the sense of Borel-Casselman \cite{BC}. We consider also
\begin{equation}
H^\bullet (S_K) = H^\bullet (\gog, K_\infty ;\Aa(G_K))
\end{equation}
the Betti cohomology of $S_K$; here $\Aa(G_K)$ is the space of automorphic forms; and we have a natural map $\mathbb{H}^\bullet(S_K) \rg H^\bullet (S_K)$, and by definition we set
\begin{equation}
H_{(2)}^\bullet(S_K)= \mathrm{Im}(\mathbb{H}^\bullet(S_K) \lgr H^\bullet(S_K)).
\end{equation}

So far all our cohomology spaces were taken with complex coefficients. The spaces $\mathbb{H}^\bullet$ and $H^\bullet$ have a topological interpretation and are therefore defined over $\Q$. For $\mathbb{H}^\bullet$, this follows from the proof of Zucker's conjecture:  see Zucker \cite{Z} for an exposition. We have $\mathbb{H}^\bullet(S_K)=IH^\bullet(S_K^*)$, a natural isomorphism, where $S_K^*$ is the Baily-Borel compactification; the (complex) intersection cohomology $\IH^\bullet$ (with middle perversity) has a natural $\Q$-structure. We now denote by $\IH^\bullet_{\Q}$,  $H^\bullet_ {\Q}$ the cohomology spaces with $\Q$-coefficients; and $\mathbb{H}^\bullet_{\Q}$ the $\Q$-space defined  through the isomorphism with $\IH^\bullet$.


Let $\Hh_K = \Hh$ be the Hecke algebra of level $K$ (with $\Q$-coefficients). It acts naturally on $\mathbb{H}^{\bullet}(S_K)$.

\vspace{2mm}

{\maprop
(i) $\Hh_K $ acts naturally on $\IH^{\bullet}_{\Q}(S_K^*)$.

(ii) The isomorphism $\mathbb{H}^{\bullet}(S_K) \rightarrow \IH^{\bullet}(S_K^*)$ is Hecke-equivariant.

(iii) There is a natural map $J: \IH^{\bullet}_{\Q}(S_K^*) \rightarrow H^{\bullet}_{\Q} (S_K)$ which, tensored with $\C$, induces the natural map $\mathbb{H}^{\bullet} (S_K) \rightarrow H^{\bullet}(S_K, \C).$}

\vspace{2mm}

Proof: Ad(i,ii): See Loojenga\cite [p.18]{Loo}. Note that Loojenga  gives his proof for classical Hecke correspondences associated to double cosets in $G(\Q)$ with respect to a congruence subgroup. Since $G$ is semisimple and simply connected, these coincide with the ad\`elic Hecke correspondences.

Ad(iii):  The intersection cohomology is the cohomology of a complex of sheaves $IC_{S*}$ of $\Q$-vector spaces \cite{BIC}. If $j: S \rightarrow S^*$ (we drop the index $K$ for simplicity) then $j^* IC_{S^*}$ is quasi-isomorphic to $\Q_S$. Therefore we obtain a natural map $\IH^{\bullet}_{\Q}(S^*) = H^{\bullet}(S^*, IC_{S^*}) \rightarrow H^{\bullet}(S, \Q)$. On the other hand, Zucker defined a complex of fine (complex) sheaves $\Ll_{(2)}$ on $S^*$. The basic result is that $\Ll_{(2)}$  and $IC_{S^*} \otimes \C$ are quasi-isomorphic. Thus $\IH^{\bullet} (S_K)$ is defined by the complex of global sections of $\Ll_{(2)}$ over $S^*$, and this is the cohomology of the space of global $L^2$-forms, i.e. $\mathbb{H}^{\bullet}(S)$. Over $\C$, the natural map  $\mathbb{H}^{\bullet}(S) \rightarrow H^{\bullet} (S, \C)$  is then defined by the restriction of differential forms.

\vspace{2mm}

From (ii) we deduce that $\Hh_K$ preserves $\mathbb{H}^{\bullet}_{\Q}(S)$.
If we define $\mathbb{H}^\bullet(\tilde{S})$ as $\liminj\limits_K \mathbb{H}^\bullet(S_K)$ we see that $\Ll(\pi_f,1/2)=\bigotimes\limits_{v\nmid\infty} \Ll(\pi_r,\frac{1}{2})$ is a $G(\A_f)$-submodule of $\mathbb{H}^q(\tilde{S})$ for $q=\sum q_v(v|\infty)$. Probably it occurs with multiplicity $1$ (see below). Equivalently, if $\pi_v$ is unramified for $v\notin S$ and $K_v(v\in S)$ is a suitable subgroup having a triangular decomposition, the representation of $\bigotimes\limits_v \Hh_v(G_v,K_v)$ occurs in $\mathbb{H}^\bullet(S_K)$.

The image of $\Ll(\pi_f,1/2)^K$ in $\mathbb{H}^\bullet(S_K)$ has a further property. Let $S_K^{BS}$ be the Borel-Serre compactification of $S_K$ (\cite{BS}; see \cite{CE} for the ad\`{e}lic  formulation.) It contains in particular a facet $e'(P)$, with $P=MN$ and $e'(P)=A_M\ba P (\A)/K_P\, K_{M,\infty}$ $K_P\subset P(\A_f)$ compact-open, and $M\subset P(\R)$ a Levi subgroup invariant by the Cartan involution associated to $K_\infty$. Cf.~\cite{CE} \footnote{In fact $e'(P)$ is a disjoint union ot the facets denoted by $e'(P)$ in \cite{BS}} . Here $A_M=S_{M_0}(\R)^+$ where $M_0= \Res_{F/ \Q} M$ and  $S_{M_0}$ is the split centre of~$M$.

We have $H^\bullet(S_K,\C)=H^\bullet(S_K^{BS},\C)$, thus we can consider the restriction
$$
\Res:H^\bullet (S_K,\C) \lgr H^\bullet(e'(P),\C).
$$
It is of course defined over $\Q$. Now we have: \cite[Theorem III.1] {RS}:

\vspace{2mm}

{\montheo (Rohlfs, Speh) If $\omega$ is a form in $\mathbb{H}^q(S_K)$ realised as a residue from $(\pi,1/2)$, $\Res\circ J(\omega)\not=0$.}

\vspace{2mm}

See \cite{RS}. Now let $a\in \Aut(\C/\Q)$. Then, as recalled in Chapter~1, $a(\pi)$ is naturally defined. It has the same properties as $\pi$. In particular Theorem~2.3 applies to it, $\psi$ (\S~2.1) being now replaced by
$$
a(\psi)=a(\pi) \otimes \usp(2)\oplus \bf{1}.
$$

Let $\pi_{G,f}=\Ll(\pi_f,1/2) \subset \mathbb{H}^q(\tilde{S})$. We can consider its transform by $a$. Since $\mathbb{H}^q(\tilde{S})$ carries a $\Q$-structure conserved by $G(\A_f)$,  this transform is contained in $\mathbb{H}^q(\tilde{S})$ .

\vspace{2mm}

{\monlem $a(\pi_{G,f}) \cong\Ll(a(\pi_f),1/2)$.}

\vspace{2mm}

This follows easily from the methods in \cite[\S~3.2]{Cl}. For finite $v$, $\pi_{G,v}$ is the unique irreducible quotient of $\ind (\pi_v |\,|^{1/2})$. The induction is \textit{unitary} induction and is not a rational operation. However $\ind (?)=\Ind(? \delta_ P^{1/2})$ where $\Ind$ denotes non-normalised induction and $\delta_P^{1/2}(m) = |\det m|^{\frac{n+1}{2}}$. Thus $\ind(\pi_v |\,|^{1/2}) = \Ind(\pi_v |\,|^{\frac{n}{2}+1})$. Since $n$ is even, we see that $\ind(a( \pi_v)|\,|^{1/2})=a(\ind(\pi_v|\,|^{1/2}))$. Since $\Ll$ is the unique irreducible quotient, the same applies to~it.

Since $\mathbb{H}^q$ carries a $\Q$-structure conserved by $G(\A_f)$, we see that $\Ll(a(\pi_f),1/2)$ is contained in $\mathbb{H}^q(\tilde{S})$; it comes from a representation $a(\pi_G)$ in the discrete spectrum, associated to $a(\psi)$. Moreover $a(\pi_G)$ is not cuspidal: if $\omega$ belongs to our given copy of $\Ll(a(\pi_f),1/2)$, $\Res \circ J(\omega)\neq 0$ and, $\Res \circ J$ being defined over $\Q$, $\Res \circ J(a(\omega))\neq 0$. Consequently there exists a \textit{residual} representation in $\Hh_{a(\psi)}$. By the previous considerations,  this implies that $L(1/2,a(\pi)) \neq 0$.

Note that $L(s,a(\pi),\Lambda^2)$ has a pole at $s=1$ by Proposition~1.1. However this is in fact implied by the foregoing argument. We note that the sign $\varepsilon(1/2,\pi)$ (cf.  Proposition~2.1) has not figured in the argument. This is not surprising, because of one of the properties of Arthur's construction. Indeed we had a representation $\pi_G$ of the form $\Ll(\pi_v,1/2)$ at all primes. For this Langlands quotient the sign $\langle \cdot,\pi_{G,v}\rangle$ (cf.~\S~2.1) is trivial for all $v$. This follows from Arthur's  \cite[Prop. 7.4.1]{A}  and the fact that the $L$-packet of $\pi$ is reduced to one element --- the group being $\GL(n)$. Thus $\Ll(\pi,1/2)$ does not occur if $\varepsilon(1/2,\pi)=-1$. In any case we have
$$
\begin{array}{ccc}
L(1/2,\pi)\neq 0 &\Leftrightarrow &L(1/2,a(\pi)) \neq 0\\
\Downarrow &&\Downarrow\\
\varepsilon(1/2,\pi)=1 &&\varepsilon(1/2,a(\pi))=1.
\end{array}
$$

This leads one to conjecture (but does not prove):

\vspace{2mm}

{\maconj ($\pi$ superregular algebraic, self-dual representation of\break $\GL(2n,\A_F)$.) For $a\in \Aut(\C/\Q)$, }
$$
\varepsilon(1/2,\pi) = \varepsilon(1/2,a(\pi)).
$$

We  have now proven (at least when $\pi_G$ has cohomology with trivial coefficients):

\vspace{2mm}

{\montheo Assume $\pi$ a representation of $\GL(2n,\A_f)$, is superregular algebraic and self-dual. Then, for any $a\in \Aut(\C/\Q)$, $L(1/2,\pi)=0\Leftrightarrow L(1/2,a(\pi))=0$.}

\vspace{2mm}

\subsection{}

We now explain how to extend this argument to the case where $\pi_G$ has cohomology with non-trivial coefficients.

(a) Let $G_0= \Res_{F/ \Q}G$. Then $G_{\infty} = G_0(\R) = \prod_{\iota} G(F_{\iota}\otimes \R)$ where $\iota$ runs over the embeddings $F \rightarrow  \R$ and $F_{\iota}= \iota(F) \subset \bar{\Q}$. We first consider a complex coefficient system, i.e. an irreducible representation $V$ of $G_0(\R)$. Thus $V= \bigotimes_{\iota} V_{\iota}$,  where $V_{\iota}$ is determined by the highest weight $m_{\iota} = m_v$ for $G(\C)$, $m_v= p' - \rho$,  where
$$
p'= (p_1+1/2, p_1-1/2, \dots, p_m+1/2, p_m-1/2), 
$$
 cf. the formula before Proposition 2.1, and $\rho = (n, n-1, \dots , 1)$ is the half-sum of roots of $G$ (for the implicit choice of ordering.) Here $(p_i)$ is defined of course by $v$ or, equivalently, $\iota$; $V=\bigotimes V_{\iota}^{m_{\iota}}$ with obvious notation.
 
 \vspace{2mm}
 
 (b) Now $G_0(\Q)= G(F) \subset \prod G(F_{\iota}).$ An irreducible representation of $G$ is always defined over $\Q$, cf. Serre \cite{Se2}, so $G(F)$ operates on $V_{\iota}^{m_{\iota}}$ by endomorphisms defined over $\iota(F) \subset \C$. Thus the action of $G(F)$ on $V$ is defined over $\prod \iota(F) = F_{gal} \subset \C$, and the same is true of the local system $\mathcal{V}$ defined by $V$.
 
  \vspace{2mm}
  
  (c) We now have an intersection cohomology sheaf $IC_{S^*}(\mathcal{V})$ with coefficients in $F_{gal}$. (We simply write $S$ for $S_K$.) The maps used in the previous argument,
  $$
  \IH^{\bullet}(S^*, \mathcal{V}) \rightarrow H^{\bullet}(S, \mathcal{V}) \rightarrow H^{\bullet}(e'(P),  \mathcal{V}) 
  $$
are all maps of $F_{gal}$-vector spaces. Finally, $\IH^{\bullet}(S^*, \mathcal{V})$ defines, through Zucker's isomorphism, an $F_{gal}$-structure on the $L^2$-cohomology  $\mathbb{H}^q(S_K, \mathcal{V})$.
    
 \vspace{2mm}
  
  (d) However, we want to obtain a result concerning $a(\pi_{G,f})$ for $a \in \Aut(\C /  \Q)$, not only in $\Aut (\C / F_{gal}).$  The action of  $\Aut(\C /  \Q)$  on $\pi_{G, \infty}$ is determined by the 'reciprocity law at infinity', cf. \cite [Theorem 3.13]{Cl} in the case of $\GL(n)$. In fact, an automorphism $a \in \Aut(\C / \Q)$ transforms the representation $V = \bigotimes V^{m_{\iota}}$ into $a(V) = \bigotimes V^{m(a \iota)} $. The associated representation  $a(\pi_{G,f})$ will therefore come from a global representation having a type at infinity permuted from that of $\pi_G$ by the permutation $\iota \rightarrow a \iota.$
  
  Therefore, consider the reducible coefficient system $\mathcal{W}= \bigoplus \sigma (\mathcal{V})$, where $\sigma \in \Gal(F/ \Q)$. Now $\mathcal{W}$ is defined over $\Q$.  The remarks in (c) apply, the local system, and $IC_{S*}(\mathcal{W})$, being defined over $\Q$. We then deduce that  $a(\pi_{G,f})$ will occur in $\mathbb{H}^q(S_K, \mathcal{W}).$ The argument concerning the restriction to $e'(P)$ remains correct, so $a(\pi_{G,f})$ will originate from a residual representation. We are now reduced to the foregoing proof.
  
   \vspace{2mm}

We recall again that Theorem 2.5 follows from the explicit formula of Grobner and Raghuram. The present proof (for this special result) is more direct. However, we will now consider other cases.

\section{Rankin $L$-functions}

\subsection{}
In this chapter we extend the method, and the main result, of Chapter~2, to Rankin $L$--functions. We consider  $RASDC$ representations $\pi$, $\rho$, respectively, of $\GL(n_1,\A)$ and $\GL(n_2,\A)$, and the Rankin $L$--function
\begin{equation}
L(s,\pi \times \rho).
\end{equation}
As in Chapter 1 this is the complete $L$-function. We can consider $a(\pi)$, $a(\rho)$ for $a\in \Aut(\C/\Q)$ and are interested in the ``invariance of the zero'' at~$s=1/2$, assuming $n_1 n_2$ even. 

We want to consider Eisenstein series (for a maximal parabolic subgroup) in the variable $s$, such that the  numerator of the constant term is essentially
$$L(s,\pi\times\rho) L(2s,\pi,R)
$$
where $L(2s,\pi,R)$ should have a pole at $s=1$, and $R$, a representation of the dual group $\GL(n_1,\C)$, will be $\Lambda^2$ or~$\mathrm{Sym}^2$.

Since we are relying on Arthur's results the ``ambient'' group $G$ must be split classical. Furthermore we must use the crucial rationality result (Proposition~2.3) and we need Zucker's conjecture or its generalisation by Saper \cite{Sap}. In particular, following Langlands \cite{LEP} we see that there are two possible cases. We refer to Kim \cite{HK} for a useful exposition of Langlands' results.

\subsection{} In the first case we consider, as before, $G = \Sp(n)$, but the maximal parabolic subgroup has Levi subgroup $M = \GL(r)\times \Sp(n-r)$. We are given a cuspidal representation $\pi\otimes \sigma$ of~$M(\A$).

We have $\hat{G}= \SO(2n+1,\C)$, containing 
$$
\hat{M} = \GL(r,\C) \times \SO(2n+1-2r,\C).
$$
Setting $N'=2(n-r)+1$, we will assume that $\sigma$ is (by Arthur's theorems) associated to a cuspidal, self-dual representation $\rho$ of $\GL(N',\A)$. The constant term of the Eisenstein series constructed from
$$
\ind _{MN(\A)}^{G(\A)} \pi[s] \otimes \sigma
$$
is given, as usual, by a ratio with numerator
\begin{equation}
L(s,\pi\times \rho)L(2s,\pi,\Lambda^2).
\end{equation}
Cf.~\cite[p.143]{HK}. (If $r=1$, only $L(s,\pi\times \rho)$ occurs.). In particular the interesting case is that where $\pi$ is symplectic. Thus we have
\begin{equation}
L(s,\pi \times \rho),
\end{equation}
\begin{equation}
\left\{
\begin{array}{ll}
\pi\ RASDC\ \mathrm{representation\ of\ }& \GL(r,\A),\\
&r\ \mathrm{even\ }\ge 2\\
\rho\ RASDC\ \mathrm{representation\ of\ } & \GL(t,\A),\\
&t\ \mathrm{odd}.
\end{array}
\right.
\end{equation}

\subsection{} 
The other interesting case occurs for the split group $G=\SO(2n+1)$. For any real prime $v$, $G\times F_v$ is of equal rank, i.e. $\rk_\R G= \rk_\R \,K_\infty$. The extension of Zucker's conjecture has been proven by Saper. Again we take $M = \GL(r) \times \SO(2(n-r)+1)$, and $\pi \otimes \sigma$ an irreducible cuspidal representation.

We have $\hat{G} = \Sp(n,\C)$, containing $\hat{M} = \GL(r,\C) \times \Sp(n-r,\C)$. We assume that $\sigma$ comes from a cuspidal, self-dual representation $\rho$ of $\GL(N',\A)$, $N'=2(n-r)$. Then $\rho$ has to be of symplectic type. The relevant numerator is then, for $r<n$
\begin{equation}
L(s,\pi\times\rho)L(2s,\pi,S^2).
\end{equation}

Thus $\pi$ has to be orthogonal. We see that this will yield no new results compared with (3.4). Also note that (3.5) is incorrect if $r=n$, in which case the numerator~is
$$
L(2s,\pi,S^2).
$$

This could be used to show that ``$\pi$ orthogonal'' is invariant by $\Aut(\C/\Q)$, but we have seen that this is known for regular algebraic representations.

In conclusion, we will only consider Eisenstein series for $\Sp(n)$. We have recorded this calculation since other consequences might follow from this construction.

\subsection{} From now on our datum will be a  pair $\pi$, $\rho$ verifying (3.4). We set $t=2m+1$, $n=r+m$; $\pi$ is symplectic, $\rho$ is orthogonal. Also set $r=2r_1$.

Before we proceed it will be necessary to review the correspondence between $\sigma$, a cuspidal representation of $\Sp(n-r,\A)$, and $\rho$, a cuspidal, self-dual, orthogonal representation of $\GL(N',\A)$, $N'=2(n-r)+1$. We have \textit{assumed} $\rho$ cuspidal. Thus $\rho\in \tilde{\Psi}_{sim}(G')$ where $G'=\Sp(n-r)$. Cf.~Arthur \cite[p.~37]{A}. By Arthur's Theorem~1.5.2, the parameter $\psi'\equiv \rho$ defines a collection of representation $\sigma$ occurring in $L_{disc}^2(G'(F)\ba G'(\A))$.

\vspace{2mm}

{\maprop Any $\sigma$ in the A-packet defined by $\rho$ is cuspidal. Moreover this A-packet is non-empty.}

\vspace{2mm}

Indeed we know that, $\rho$ being $RASCD$, it is tempered at all primes. The factors $\sigma_v$ of $\sigma$ are therefore tempered at all primes, by \cite[Thm.1.5.1]{A}. Since $\sigma$ occurs in the discrete spectrum, it is cuspidal by a well-known result.

Moreover, in this case, $\Ll_{\psi'} = \SO(N') \subset \GL(N')$ \cite [1.4.4]{A}, and $\CSs_{\psi'} = S_{\psi'}=\{1\}$. This implies that a tensor product $\sigma=\otimes \sigma_v$, with $\sigma_v$ any local representation in the $\psi'$-packet, occurs in $L_{disc}^2 (G(F)\ba G(\A))_{\psi'}$.

\textbf{Remark.} (1) By construction the first factor of (3.2) (coming from the theory of Eisenstein series) considered outside a finite number $S$ of primes, is equal to $L^S(s,\pi\times\rho)$.

(2) The case of $G'=\SL(2)$ shows of  that the set $\{\sigma\}$ is, in general, infinite.

\subsection{} We are now interested in the $\psi$-parameter (for~$G$)
$$
\psi=\pi \otimes \usp(2) \oplus \rho.
$$
Since $\pi$ is symplectic, as well as $\usp(2)$, this is a parameter associated to~$G$.

As in Chapter 2, it defines (almost everywhere) Hecke matrices (in $\hat{G}$) which coincide with those of the (ad\`elic) induced representation
\begin{equation}
\ind _{MN}^G (\pi[s]\otimes \sigma)\qquad (s=1/2).
\end{equation}
We first have to determine when such a representation will be cohomological. We consider a place~$v|\infty$.

The infinitesimal character if $\pi_v$ is of the form $p=(p_1,\ldots p_r)$, with
$$
p_1>p_2\cdots>p_{r_1}>-p_{r_1}> \cdots> -p_1
$$
with $p_i\in\frac{1}{2}+\Z$. For $s=1/2$, we see that the infinitesimal character of $\pi[1/2]$~is
$$
Y_1=\{p_1+1/2, \ldots, p_{r_1}+1/2,-p_{r_1}+1/2,\ldots -p_1+1/2\}.
$$
It contributes a term $(Y_1,-Y_1)$ to the character of the representation (3.6) lifted to $\GL(N)$. (See~\S~2.2.)

On the other hand, the infinitesimal character of $\sigma$ can be seen, lifted to $\hat{G}'$, as the infinitesimal character of $\rho$, of the form
$$
q_1>q_2>\cdots q_m >0>-q_m>\cdots -q_1\,.
$$
Finally, the infinitesimal character of the induced representation (3.6) is obtained as a parameter of dimension $2n+1$, equal~to
\begin{equation}
\{\pm p_i \pm 1/2,\ \pm q_j,\, 0\},
\end{equation}
and defined modulo $\Ss_N$, $N=2n+1$. It must be regular, so we must assume
$$
\left\{
\begin{array}{l}
p_i > p_{i+1}+1,\ p_{r_1}\geq \frac{3}{2}\\ \\
p_i \pm \frac{1}{2}\not= q_j\quad (\mathrm{all\ }i,j).
\end{array}
\right.
$$
Modulo $\Ss_N$, we see that the character is $(Y,0,-Y)$ where
$$
\begin{array}{l}
Y=(p_1+1/2,p_1-1/2,\ldots, q_1, \ldots,p_{a_1}+1/2,p_{a_2}+1/2,\ldots \dots , q_2,\ldots\\
\vspace{2mm}
\hspace{7cm}\ldots , q_m,p_{r_1}+1/2, p_{r_1}-1/2)
\end{array} 
$$
the $q_j$ being interspersed so that
$$
Y=(m_1,\ldots m_n)
$$
with $m_1>\cdots > m_n >0$.  Here, $Y$ is the infinitesimal character of $\ind (\pi[1/2]\otimes \sigma)$.

\vspace{2mm}

{\madef Assume $\pi$ is superregular,  and $\rho$ is regular. We say that $(\pi ,\rho)$ are \emph{disjoint} if, for all $i=1,\ldots r_1$ and $j=1,\ldots m$,
$$
p_i\pm 1/2 \neq q_j.
$$

In this case the infinitesimal character of $\ind (\pi[1/2]\otimes \sigma)$ is regular.}

\vspace{2mm}

\textbf{Remark.}  We will be considering the Rankin $L$-function $L(s,\pi\times \rho)$. Conjecturally this is the standard $L$-function $L(s,\pi \boxtimes\rho)$ where $\pi \boxtimes \rho$ is the Langlands tensor product of $\pi$, $\rho$. This representation should be symplectic. Recalling that the parameters $p,q$ are associated to the representations $z\mapsto(z/\bar{z})^p$, $z\mapsto (z/\bar{z})^q$ of $W_\C$, we see that the Langlands parameters (or the infinitesimal parameters) of $\pi \boxtimes \rho$ will be, with obvious notation,
$$
(p_1,\ldots -p_1) \otimes (q_1,\ldots 0,\cdots -q_1).
$$
In particular $\pi\boxtimes \rho$ will be regular if

\vspace{2mm}

(i) the $\pm p_i \pm q_j$ are distinct

(ii) $\pm p_i \neq \pm p_i \pm q_j$.

\vspace{2mm}

This is a much stronger condition. For instance, assume $m=1$, $r_1=2$, and we want to ensure the condition (with~$q>0$)
$$
p_1+1/2 > p_1-1/2 > q >p_2 +1/2 > p_2 -1/2  \geq 1.
$$
Can the tensor product be singular, e.g. $p_1-q=p_2+q$~?  We set $p_1=p_2+2q$ and we must ensure
$$
p_2+2q+1/2 > p_2 +2q-1/2 > q > p_2+1/2 > p_2 -1/2 \ge 1.
$$
Thus we seek $p_2 >1/2 -q$, $p_2<q-1/2$ (and $p_2 \ge 3/2$). This is satisfied for $q$ large.

\subsection{} In this excursus, we return to the relation between $\rho$ and $\sigma$ (cf. Proposition~3.1) and consider questions of rationality. We first make some general remarks, in the case where (the generalisation of) Zucker's conjecture is true. Changing notation, we denote by $G$ a semi-simple group over a totally real field $F$, such that for each $v|\infty$ $\rk(G_v)=\rk(K_v)$ where $K_v\subset G(F_v)$ is maximal compact. If $G$ is Hermitian symmetric, or if $G$ is of type $B_n$, there is an isomorphism
$$
\mathbb{H}^\bullet(S_K) \cong \IH^\bullet (S_K^*)
$$
where $S_K$ is defined as in (2.7) --- $K\subset G(\A_f)$ is a compact subgroup --- and $\IH^\bullet(S_K^*)$ is the middle perversity intersection cohomology  of a suitable compactification of $S_K$. See \cite{Sap} for the non--Hermitian case.

We deduce from this a natural $\Q$--structure on $\mathbb{H}^\bullet(S_K)$ or $\mathbb{H}^\bullet(\tilde{S}) =\break \liminj\limits_K \mathbb{H}^\bullet (S_K)$. We have
$$
\mathbb{H}^\bullet(\tilde{S}) = \bigoplus_\pi H^{\bullet} (\gog,K_{\infty};\pi_\infty) \otimes \pi_f
$$
where $\bigoplus \pi$ is the discrete spectrum. In  particular (in all these  cases) this direct sum is ``defined over $\Q$''. In particular, if $a\in \Aut(\C/\Q)$:
\begin{eqnarray}
&\pi_f\ \textit{representation\ of\ } G(\A_f),\ \pi_\infty \otimes \pi_f \subset L_{disc}^2(G(\Q)\ba G(\A))\nonumber\\
& \textit{for\ some\ cohomological\ } \pi_\infty \Rightarrow \\
&\textit{$a(\pi_f)$ has\ the\ same\ property\ for\ any}\ ~~a\in \Aut(\C/\Q).\nonumber
\end{eqnarray}
Here $a\in \Aut(\C)$ acts naturally, as in Chapter~1. (We should point out, however, that if $G$ is not $\GL(n)$, the field of rationality and (a) field of definition may not coincide --- for $\GL(n)$ see \cite[\S~3.1]{Cl}. 

Note too that if $\pi=\pi_\infty \otimes \pi_f$ is cuspidal, this does not imply that $\pi'=\pi_\infty'\otimes a(\pi_f)$ is~cuspidal. 

\subsection{} We now return to our initial situation and apply this to the representation $\sigma$ of $\Sp(m,\A)$. We have seen that it was cuspidal (Proposition~3.1). For $a\in \mathrm{Aut(\C})$, consider a representation $\sigma'=\sigma_\infty'\otimes a(\sigma_f)$ so obtained; we can assume
$$
H^q(\gog, K_\infty;\sigma_\infty') \neq 0
$$
if $H^q(\gog,K_\infty;\sigma_\infty) \neq 0$.

The representation $\sigma$ belongs to the Arthur-packet defined by $\rho$. The Arthur packets are disjoint (\cite[Thm.~1.5.2]{A}\footnote{The global packets, not their local components!} and defined by their Hecke matrices at almost all primes. (This is easily seen by an adaptation of the proof of~Thm.~2.3.).

Let $\Hh_\rho$ denote, as in \S~2.1, the subspace of $L_{disc}^2(G'(F)\ba G'(\A))$ contributed by $\psi'=\rho$.

\vspace{2mm}

{\monlem $\sigma_\infty \otimes \sigma_f\subset \Hh_\rho$ (for some $\sigma_\infty$) if, and only if, $\sigma_\infty' \otimes a(\sigma_f)\subset \Hh_{a(\rho)}$ for some $\sigma_\infty'$.}

\vspace{2mm}

Let $S$ be a large set of primes containing $\infty$. Then $\sigma_\infty \otimes \sigma_f\subset \Hh_\rho$ if, and only if, for~$v\notin S$:
\begin{equation}
t_{\rho,v} = t_{\sigma,v} \oplus 1\oplus t_{\tilde{\sigma},v}
\end{equation}
with obvious notation. These are the Hecke matrices for the Langlands normalisation.

Here, for instance, $t_{\sigma,v} \in \hat{T}_{G'}$ and $\hat{T}_{G'}\cong \G_m(\C)^n$ has an obvious $\Q$-structure, as does $\hat{T}_{\GL(N)}$.

Define $t_{\rho,v}^T = t_{\rho,v}q_v^{\frac{1-N'}{2}}$ (Tate normalisation: cf. \cite[p.~106--108]{Cl}).

We then have $t_{a(\rho),v}^T= a(t_{\rho,v}
^T)$, cf.~\cite[Lemma~3.5]{Cl}. Since $N'$ is odd, this implies that $t_{a(\rho),v}=a(t_{\rho,v})$.

On the other hand, $t_{\sigma,v}$ is defined by
$$
\varphi\cdot e_v = \widehat{S \varphi}(t_{\sigma,v})
$$
where $e_v$, in the space of  $ \sigma_v$, is the unramified vector, $\varphi$ belongs to the local unramified Hecke algebra for $\Sp(m)$, $S\varphi$ is the function on $T_{G'}(F_v)$ given~by
$$
S\varphi(\tau) = \delta_B^{1/2}(\tau) \, \int_{N(F_v)} f(\tau n) dn
$$
and $\delta_B^{1/2}$, $B=T_{G'}N_{G'}$, is the half modulus. \cite[p.~146]{Car}. Finally, $S\varphi=\sum\limits_{\lambda\in X_*(T_{G'})} b_\lambda ch(\varpi^\lambda T_{G'}( \Oo))$  and $\widehat{S \varphi}(t)=\sum b_\lambda t^\lambda$.  Since for $\Sp(m)$ $\delta_B^{1/2}$ is a rational character, this implies
$$
S(a(\varphi)) =a(S\varphi)
$$
where $a$ acts on the coefficients; dually we deduce that $t_{a(\sigma),v}=a(t_{\sigma,v}
)$.

Finally, the relation (3.9) is preserved by the action of $a$ on the both sides, which implies the lemma.

Since $\rho$ is $RASD$ cuspidal, so is $a(\rho)$, and $a(\rho)$ satisfies the Ramanujan conjecture. We deduce from Proposition~3.1 that $a(\sigma)$ is (tempered) cuspidal.

\subsection{} 
Consider now, as in \S 3.4, a pair $(\pi,\rho)$ verifying~(3.4), and a representation $\sigma$ of $\Sp(m,\A)$ associated to $\rho$; $\sigma$ is cuspidal (Proposition 3.1). The representation $\ind _{MN}^G(\pi[1/2]\otimes \sigma)$, $M=\GL(r)\times \Sp(m)$, ad\`elic induced representation, has Hecke matrices given by $\pi \otimes \usp(2)\oplus \rho$.

Forming Eisenstein series
$$
E(-,s):\ind _{MN}^G (\pi [s]\otimes \sigma) \lgr \Aa(G(F)\ba G(\A)),
$$
absolutely convergent for $\Re(s)>\frac{n+1}{2}$ \cite[p.~86]{MW}, we see that for $f_s\in I_s =\ind (\pi[s]\otimes \sigma)$ the constant term~is
$$
E_P(f_s) = f_s + M(s)f_{-s}
$$
with $M_v(s)$ given, for $f_v$ the unramified vector in~$I_{v,s}$,~by
$$
\frac{L_v(s,\pi\times \rho) L_v(2s,\pi,\Lambda^2)}{L_v(s+1,\pi\times \rho) L_v(2s+1,\pi,\Lambda^2)}
$$

This applies for $v\notin S$, $S\supset \infty$. \textit{Assume now that} $L(1/2,\pi\times\rho)\neq 0$. Writing $M(s) = M_S(s)M^S(s)$, we see that $E(-,s)$ has a residue at $s=1/2$; we obtain an intertwining operator (cf.~\S~2.3)
$$
R(1/2) : I_{1/2} \lgr L_{disc}^2(G(F)\ba G(\A)).
$$
The representation $I_{v,1/2}$, $\pi_v$ and $\sigma_v$ being tempered at all $v$, has a unique irreducible quotient $\Ll(\pi_v,\sigma_v,1/2)$; the representation of $G(\A)$ on $L_{disc}^2$ being semi-simple, we see~that
$$
R(1/2) : \Ll(\pi,\sigma,1/2) \hookrightarrow L_{disc}^2 (G(F)\ba G(\A))
$$
with $\Ll$ the global Langlands quotient. For each Archimedean prime, $\Ll(\pi_v,\sigma_v,1/2)$ is cohomological; let $q_v$ the minimal degree of its cohomology; with $q=\sum\limits_{v|\infty}q_v$ we finally have (these cohomology spaces being one-dimensional).
\begin{equation}
\Ll_f(\pi,\sigma,1/2) \hookrightarrow H^q(\frak{g}, K_\infty; L_{disc}^2 (G(F)\ba G(\A)).
\end{equation}
The latter space is identified by Zucker's conjecture with $\IH^q(\tilde{S}^*)$. It is naturally defined over~$\Q$.

Moreover (cf. \S~2.4), if we fix a level $K\subset G(\A_f)$ and consider (3.10) at a finite~level:
$$
R(1/2):\Ll_f(\pi,\sigma,1/2)^K \hookrightarrow \IH^q(S_K^*),
$$
we know that there is a component $e'(P)$ of the boundary of the Borel--Serre compactification such that 
$\Res \circ  J(\omega) \not=0$ for some $\omega $ in the image. (The notation is that of Proposition 1.3.)

Now let $a\in \Aut(\C/\Q)$.

\vspace{2mm}

{\monlem 
$$
\Ll_f(a(\pi),a(\sigma),1/2) \cong a(\Ll_f(\pi,\sigma,1/2).
$$}

\vspace{2mm}

Note that $a(\pi)$, $a(\rho)$ are still (everywhere) tempered, and that the Langlands quotient on the left--hand side is defined in the same manner. In particular, it suffices to check that, for any finite~$v$,
$$
\ind (a(\pi),a(\sigma),1/2)=a(\ind(\pi,\sigma,1/2))
$$
with obvious notation. This would be true for non-normalised induction (without the factor $\vert ~\vert^{1/2}$). We must check that, for $m=(g,h)\in \GL(r)\times \Sp(m)$,
$$
|\det g|^{1/2} \delta_P(m)^{1/2} \ \mathrm{takes\ rational\ values,}
$$
where $\delta_P(m)=|\det(Ad(m) \vert \fn)|$. We can take $m$ in the maximal torus $T\subset M$. We have $\fn_0=\fn+\fn_M$ where $\fn_0$ is the  unipotent radical of the Lie algebra  of the Borel subgroup of $G$ and $\fn_M$ is maximal nilpotent (positive) in $\Lie(M)$,~so
$$
\delta_P(t) = \det(Ad~t\mid \fn_0):\det(Ad~t\mid \fn_M),
$$
$t\in T$. Now, for $t=(t_1,\ldots t_n)$:
$$
\begin{array}{c}
|\det (Ad(t)\mid \fn_0)|^{1/2} = |t_1^n \cdots t_n|\\ \vspace{2mm}
|\det(Ad(t)\mid \fn_M)|^{1/2}=|t_1^{\frac{r-1}{2}}\cdots t_r^{\frac{1-r}{2}}|\ |t_{r+1}^m \cdots t_n|.
\end{array}
$$
Since $r$ is even, multiplying by  $|\det g|^{1/2}=|t_1\cdots t_r|^{1/2}$ yields a rational value.

\vspace{2mm}

At this point we know that $\Ll_f(a(\pi),a(\sigma),1/2)$ occurs in $\IH^q(\tilde{S}^*)$, and moreover that it restricts non trivially to the boundary of the Borel-Serre compactification of $\tilde{S}$, cf. $\S 2.4$. (For psychological comfort consider a finite level!) Furthermore, the representation $\Ll(\pi,\sigma,1/2)$ belongs  to Arthur's $\psi$-packet with
$$
\psi=\pi \otimes \usp(2)\oplus \rho.
$$
This implies that $\Ll(a(\pi),a(\sigma),1/2)$ belongs to Arthur's packet~for
$$
a(\psi) = a(\pi)\otimes \usp(2) \oplus a(\rho).
$$
We know that it is residual since the restriction (in cohomology) to the boundary does not vanish. Renaming $a(\psi)$ as $\psi$, we see that we need only amplify the assertion of Theorem~2.3 in this context:

\vspace{2mm}

{\montheo 
Assume $\psi=\pi \otimes \usp(2)\oplus \rho$, $\pi_G\subset \Hh_\psi$ is cohomological, not cuspidal, and obtained from multi-residues and principal values of Eisenstein series from $(M,\pi_M)$, where $M\subset G$ is a Levi subgroup and $\pi_M$ cuspidal, at a point $s\in \C^{I}$. Then $M = \GL(r) \times \Sp(m)$, $\pi_M=\pi\times \sigma$ (where $\sigma$, a cuspidal representation, is associated to $\rho$) and $s=1/2$.}
\vspace{2mm}

By assumption $\pi_G$ is obtained from the formation of residues from $M=\GL(r_1)\times \cdots \times \GL(r_I) \times \Sp(m')$, and representations $\pi_i$ $(i=1,\ldots I)$ and $\sigma'$. By Arthur's results $\sigma'$ is associated to a parameter $\psi'=\bigoplus\limits_j \sigma_j \otimes \usp(m_j)$, $\sigma_j$ being cuspidal self-dual. (Moreover $\sigma'$ is cuspidal). The Hecke matrices, at almost all primes, of $\pi_G$ are given~by
\begin{equation}
T=T_v=t_\pi q^{1/2} \oplus t_\pi q^{-1/2} \oplus t_\rho
\end{equation}
where $q=q_v$. (The Hecke matrices are viewed in $\GL(N,\C)$). By Theorem~2.2 (Grba\v c-Schwermer) they are also given by
\begin{equation}
T'=T_v' = \bigoplus_{i=1}^I t_i q^{s_i} \oplus t_i^{-1} q^{-s_i} \oplus \bigoplus_{j=1}^J t_{\rho_j} \otimes \usp(m_j).
\end{equation}
Here $s_i\in \frac{1}{2}\Z$, and $t_i=t(\pi_i)$; we assume that $\omega_{\pi_i}=1$ on $\R_+^\times$, so $\pi_i$ is cuspidal unitary. (If the factor $\Sp(m')$ is trivial, the second term in the sum is $1$, with $m_j=1$.)  Since $\pi$ and $\rho$ are tempered, we deduce from $T=T'$ that $m_j\le 2$, that $|s_i|=0,1/2$; if $s_i=\pm1/2$ we can assume $s_i=1/2$ (use an associate datum). Now we obtain separate equations, using as in Chapter~2 that the eigenvalues of $t_i$, $t_{\rho_j}$ are in $[-1/2+\delta,1/2-\delta]$:
\begin{eqnarray}
t_\pi &=& \bigoplus_{s_i=1/2} t_i \oplus \bigoplus_{m_j=2} t_{\rho_j},\\
t_\rho &=&\bigoplus_{s_i=0} (t_i\oplus t_i^{-1}) \oplus \bigoplus_{m_j=1}t_{\rho_j}.
\end{eqnarray}

All representations considered are cuspidal. Consider first (3.14). By the results of Jacquet- Shalika, there must be at most one term on the right-hand side. We deduce that $\{i:s_i=0\}=\emptyset$, that $\{j:m_j=1\}$ has cardinality $1$, and that $\rho_j=\rho$.

Consider now (3.13). The first possibility is that $\#\{i:s_i=1/2\}=1$, that no $\rho_j$ intervene, and $t_\pi =t_i=t_{\pi_i}$. In this case, using the previous remark, we see that
$$
T'=t_\pi q^{1/2} \oplus t_\pi q^{-1/2} \oplus t_\rho
$$
and that $\pi_G$ is induced from $\GL(r)\times \Sp(m)$ from $(\pi,\sigma)$ as sought.

The second possibility is that $I=0$, which implies that $m=n$ and that $\pi_G$ is cuspidal. In this  case (3.12) yields 
$$
T=t_\pi q^{1/2} \oplus t_\pi q^{-1/2} \oplus t_\rho=\bigoplus t_{\rho_j}\otimes \usp(m_j)=t_\rho\oplus t_{\pi} \otimes \usp(2):
$$
the only $\rho_j$ in (3.13) is equal to $\pi$; it can of course occur (and it will occur depending on the behaviour of the $L$-functions) that the $\psi$-packet only contains  cuspidal representations.

\subsection{}  We can now apply Theorem 3.1, with $\pi,\rho$ replaced by $a(\pi)$, $a(\rho)$. Since $\Ll_f(a(\pi),a(\sigma),1/2)$ occurs in $\IH^q(\tilde{S}^*)$, with non-trivial restriction to the boundary, it is given by the cohomology of (non-cuspidal) Eisenstein series. Theorem~3.1 shows that these are obtained from Eisenstein series coming from $\ind _{MN}^G(a(\pi)[s]\times \sigma')$ (and their residues) at $s=1/2$. Here  $\sigma'$, a cuspidal representation of $\Sp(m,\A)$ is associated to $a(\rho)$.  (It coincides with $a(\sigma)$ at almost all primes.) Since the cohomology in $\mathbb{H}^q$ is given by square-integrable forms, this can only occur if these Eisenstein series present a residue at $s=1/2$. We deduce, by inverting the argument at the beginning of \S~3.8, that
$$
L(s,a(\pi) \times a(\rho)) L(2s,a(\pi),\Lambda^2)
$$
must have a pole at $(s=1/2)$, and therefore that $L(1/2,a(\pi)\times a(\rho))\neq 0$. We have therefore obtained the main result of this paper. (In \S~3.8, we have implicitly assumed that $\Ll(\pi_\infty,\sigma_\infty,1/2)$ had cohomology with trivial coefficients. The corrections in the case of a general coefficient system are dealt with in Chapter~2.)

\vspace{2mm}

{\montheo Let

$\bullet$ $\pi$ be a $RASDC$, superregular, representation of $\GL(r,\A)$ $(r$ even)

$\bullet$ $\rho$ be a $RASDC$ representation of $\GL(t,\A)$ ($t$ odd)

\noindent and assume that $(\pi,\rho)$ are disjoint (Def. $3.1$). Then, for any $a\in \Aut(\C/\Q)$, $L(1/2,\pi\times\rho)=0\Leftrightarrow L(1/2,a(\pi) \times a(\rho))=0$.}


\end{document}